\documentclass[12pt]{amsart}
\usepackage{graphicx}
\usepackage{amsmath,amsthm,amssymb,enumerate}
\usepackage{euscript,mathrsfs}
\usepackage[citecolor=blue,colorlinks=true]{hyperref}
\usepackage[left=3cm,right=3cm,top=4.5cm,bottom=4.5cm]{geometry}
\usepackage{enumitem}
\setenumerate{label={\rm (\alph{*})}}
\usepackage{color}
\catcode`\@=11 \@addtoreset{equation}{section}

\catcode`\@=12

\allowdisplaybreaks

\newtheorem{Theorem}{Theorem}[section]
\newtheorem{Proposition}[Theorem]{Proposition}
\newtheorem{Lemma}[Theorem]{Lemma}
\newtheorem{Corollary}[Theorem]{Corollary}

\theoremstyle{definition}
\newtheorem{Definition}[Theorem]{Definition}

\newtheorem{Remark}[Theorem]{Remark}

\newcommand{\bTheorem}[1]{
\begin{Theorem} \label{T#1} }
\newcommand{\eT}{\end{Theorem}}

\newcommand{\bProposition}[1]{
\begin{Proposition} \label{P#1}}
\newcommand{\eP}{\end{Proposition}}

\newcommand{\bLemma}[1]{
\begin{Lemma} \label{L#1} }
\newcommand{\eL}{\end{Lemma}}

\newcommand{\bCorollary}[1]{
\begin{Corollary} \label{C#1} }
\newcommand{\eC}{\end{Corollary}}

\newcommand{\bRemark}[1]{
\begin{Remark} \label{R#1} }
\newcommand{\eR}{\end{Remark}}

\newcommand{\bDefinition}[1]{
\begin{Definition} \label{D#1} }
\newcommand{\eD}{\end{Definition}}

\newcommand{\Q}{Q}

\newcommand{\prst}{\mathcal{P}}

\newcommand{\bfphi}{\boldsymbol{\varphi}}

\newcommand{\bFormula}[1]{
\begin{equation} \label{#1}}
\newcommand{\eF}{\end{equation}}

\newcommand{\Ov}[1]{\overline{#1}}

\newcommand{\DC}{C^\infty_c}

\newcommand{\vr}{\varrho}

\newcommand{\vt}{\vartheta}

\newcommand{\vm}{\vc{m}}

\newcommand{\vc}[1]{{\bf #1}}

\newcommand{\Div}{{\rm div}_x}
\newcommand{\Grad}{\nabla_x}

\newcommand{\dx}{\,{\rm d} {x}}

\newcommand{\dt}{\,{\rm d} t }

\newcommand{\dxdt}{\dx \ \dt}
\newcommand{\intO}[1]{\int_{\Omega} #1 \ \dx}
\newcommand{\intQ}[1]{\int_{{\Q}} #1 \ \dx}

\newcommand{\D}{{\rm d}}

\definecolor{Cgrey}{rgb}{0.85,0.85,0.85}
\definecolor{Cblue}{rgb}{0.50,0.85,0.85}
\definecolor{Cred}{rgb}{1,0,0}
\definecolor{fancy}{rgb}{0.10,0.85,0.10}

\newcommand\Cbox[2]{%
    \newbox\contentbox%
    \newbox\bkgdbox%
    \setbox\contentbox\hbox to \hsize{%
        \vtop{
            \kern\columnsep
            \hbox to \hsize{%
                \kern\columnsep%
                \advance\hsize by -2\columnsep%
                \setlength{\textwidth}{\hsize}%
                \vbox{
                    \parskip=\baselineskip
                    \parindent=0bp
                    #2
                }%
                \kern\columnsep%
            }%
            \kern\columnsep%
        }%
    }%
    \setbox\bkgdbox\vbox{
        \color{#1}
        \hrule width  \wd\contentbox %
               height \ht\contentbox %
               depth  \dp\contentbox
        \color{black}
    }%
    \wd\bkgdbox=0bp%
    \vbox{\hbox to \hsize{\box\bkgdbox\box\contentbox}}%
    \vskip\baselineskip%
}

\date{}

\begin{document}


\title[Ill posedness for full Euler]{Ill posedness for the full Euler system driven by multiplicative white noise}

\author{Elisabetta Chiodaroli}
\address[E.Chiodaroli]{Dipartimento di Matematica, Universit\` a di Pisa
Via F. Buonarroti 1/c, 56127 Pisa, Italy}
\email{elisabetta.chiodaroli@unipi.it}

\author{Eduard Feireisl}
\address[E.Feireisl]{Institute of Mathematics AS CR, \v{Z}itn\'a 25, 115 67 Praha 1, Czech Republic
and Technische Universitaet Berlin, Strasse des 17.Juni, Berlin, Germany}
\email{feireisl@math.cas.cz}
\thanks{The research of E.F. leading to these results has received funding from
the Czech Sciences Foundation (GA\v CR), Grant Agreement
18--05974S. The Institute of Mathematics of the Academy of Sciences of
the Czech Republic is supported by RVO:67985840.}

\author{Franco Flandoli}
\address[E.Flandoli]{Scuola Normale Superiore di Pisa, Piazza dei Cavalieri, Pisa, Italy}
\email{franco.flandoli@sns.it}

\begin{abstract}

We consider the Euler system describing the motion of a compressible fluid driven by a multiplicative white noise. We identify a large class of initial data for which the problem is ill posed - there exist infinitely many global in time weak solutions. The solutions are adapted to the noise and satisfy the entropy admissibility criterion.

\end{abstract}

\keywords{}

\date{\today}

\maketitle


\section{Introduction}
\label{I}

Problems in continuum fluid mechanics involving perfect (inviscid) fluids are in general \emph{ill posed} in the class of weak solutions. The adaptation of the method of convex integration, developed in the pioneering work of De Lellis and Sz\' ekelyhidi 
\cite{DelSze3} that culminated by the final proof of Onsager's conjecture for the incompressible Euler system, see Isett \cite{Ise}, 
Buckmaster et al. \cite{BuDeSzVi}, produced a number of examples of non--uniqueness even in the context of compressible fluids, see 
\cite{ChiDelKre}, \cite{ChiKre}, \cite{Fei2016}, and Markfelder, Klingenberg \cite{MarKli}, among others. In accordance with the results, the initial--value problem for the Euler system is ill posed even in the class of weak solutions satisfying various kinds of admissibility criteria as the energy and/or entropy inequality. 

There is a piece of evidence that \emph{stochastic} perturbations may provide a regularizing effect on deterministically ill--posed problems, in particular those involving transport, see e.g. 
\cite{FedFla}, \cite{FlaGubPri}, \cite{FlaMauNek}. On the other hand, as shown in \cite{BrFeHo2017B}, the isentropic Euler system driven by a general additive/multiplicative white noise is ill posed in the class of weak solutions. More specifically, there 
is a large class of initial data giving rise to infinitely many solutions defined up to a positive stopping time. These solutions, however, may experience an initial energy jump and as such can be discarded as physically irrelevant. 

Our goal in the present paper 
is to show the existence of infinitely many global--in--time weak solutions to a stochastically driven Euler system that are physically admissible -- they conserve the total energy and satisfy the differential version of the entropy inequality. Specifically, we consider the \emph{full Euler system}:     
\begin{equation} \label{I1}
\begin{split}
\D \vr + \Div \vm \dt &= 0\\
\D \vm + \Div \left( \frac{\vm \otimes \vm}{\vr} \right) \dt + \Grad p \dt &= - \frac{1}{2} \vm \circ \D W\\
\D E + \Div \left( (E + p) \frac{\vm}{\vr} \right) \dt &= - E \circ \D W,
\end{split}
\end{equation}
where
\[
E = \frac{1}{2} \frac{|\vm|^2}{\vr} + \vr e, \ p = (\gamma - 1) \vr e, \ \gamma > 1.
\]
Introducing the temperature $\vt$ via Boyle--Mariotte law,
\begin{equation} \label{I2}
p = \vr \vt,\ e = \frac{1}{\gamma - 1} \vt \equiv c_v \vt,
\end{equation} 
we obtain the entropy inequality 
\begin{equation} \label{I3}
\D (\vr s) + \Div (s \vm) \dt \geq - c_v \vr \circ \D W.
\end{equation}
For definiteness, we impose the impermeability condition 
\begin{equation} \label{I4}
\vm \cdot \vc{n}|_{\partial Q} = 0.
\end{equation} 
Here,
$W$ denotes the standard scalar valued Wiener process while the symbol $\circ$ indicates that the stochastic integral in the weak formulation of the problem is interpreted in the Stratonovich sense.

We show the existence of infinitely many solutions following the strategy of Luo, Xie, Xin \cite{LuXiXi} used also in \cite{FeKlKrMa}. 
Specifically, we choose arbitrary \emph{piece--wise} constant initial distributions of the density and the absolute temperature and 
we transform the problem into a family of partial differential equations with random parameters. Then we  
apply the result of De Lellis and Sz\' ekelyhidi \cite{DelSze3} for the incompressible Euler system with constant pressure 
on each domain where the initial data are constant. Finally, 
we pass back to the original system ``pasting'' together the solutions previously obtained. The issue of progressive measurability 
of the oscillatory solutions, that was absolutely crucial for the analysis in \cite{BrFeHo2017B}, is handled here by introducing 
a new stochastically rescaled time variable. 
  
The paper is organized as follows. In Section \ref{M}, we introduce the necessary preliminary material and state our main result. 
The main ideas of the proof are described in Section \ref{H}, where the transformation into a system with random coefficients is performed. In Section \ref{C}, we apply the nowadays standard tools of convex integration to the transformed problem. In Section \ref{R}, we introduce a new ``random'' time variable. The existence proof is completed in Section \ref{S}.  	

\section{Main result}
\label{M}

Let $\{ \Omega, \mathfrak{F}, \mathcal{P} \}$ be a probability basis, with a right continuous complete filtration 
$\{ \mathfrak{F}_t \}_{t \geq 0}$, and a Wiener process $W$. 

\begin{Definition} \label{D1}

We say that $[\vr, \vm, E]$ is a weak solution of the Euler system \eqref{I1}, with the boundary condition \eqref{I4}, and the initial condition 
\begin{equation} \label{M1}
\vr(0, \cdot) = \vr_0, \ \vm(0, \cdot) = \vm_0, \ E(0, \cdot) = E_0,
\end{equation}
if:

\begin{itemize}
\item $\vr \geq 0$ $\mathcal{P}-$a.s., the functions 
\[
t \mapsto \intQ{ \vr(t, \cdot) \varphi}, \ t \mapsto \intQ{ \vm(t, \cdot) \cdot \bfphi}
\]
are continuous $\{ \mathfrak{F}_t \}_{t \geq 0}$--adapted semimartingales for any $\varphi \in C^1(\Ov{Q})$, 
$\bfphi \in C^1(\Ov{Q}; R^N)$, respectively,
\begin{equation} \label{M2}
\intQ{ \vr (\tau, \cdot) \varphi } = \intQ{ \vr_0 \varphi } + \int_0^\tau \intQ{ \vm \cdot \Grad \varphi }  
\end{equation} 
for any $\tau \geq 0$ and any $\varphi \in C^1(\overline{Q})$;
\item $E - \frac{1}{2} \frac{|\vm|^2}{\vr} \geq 0$ $\mathcal{P}$-a.s., 
the function 
\[
t \mapsto \intQ{ E(t, \cdot) \varphi }
\]
is a continuous $\{ \mathfrak{F}_t \}_{t \geq 0}$--adapted semimartingale for any $\varphi \in C^1(\Ov{Q})$;
\begin{equation} \label{M3}
\begin{split}
&\intQ{ \vm (\tau, \cdot) \cdot \bfphi } \\&= \intQ{ \vm_0 \cdot \bfphi } + \int_0^\tau \intQ{ \left[ \frac{\vm  \otimes \vm}{\vr} : \Grad \bfphi + p \Div \bfphi \right]} - \frac{1}{2} \int_0^\tau \left( \intQ{ \vm \cdot \bfphi } \right) \circ \D W   
\end{split}
\end{equation}
for any $\tau \geq 0$ and any $\bfphi \in C^1(\Ov{\Q}; R^N)$, $\bfphi \cdot \vc{n}|_{\partial Q} = 0$,
where $p = \left({\gamma - 1}\right) \left[ E - \frac{1}{2} \frac{|\vm|^2}{\vr} \right]$;
\item 
the energy equality 
\begin{equation} \label{M4}
\intQ{ E(\tau, \cdot) \varphi } = 
\intQ{ E_0 \varphi } + \int_0^\tau \intQ{ \left( E + p \right) \frac{\vm}{\vr} \cdot \Grad \varphi } \dt - 
\int_0^\tau \left( \intQ{ E \varphi } \right) \circ \D W 
\end{equation}
holds for any $\tau \geq 0$ and any $\varphi \in C^1(\Ov{Q})$;
\item
the entropy inequality 
\begin{equation} \label{M5}
\intQ{ \vr s (\tau, \cdot) \varphi } \geq \intQ{ \vr_0 s(\vr_0, E_0) \varphi }
+ \int_0^\tau \intQ{ s \vc{m} \cdot \Grad \varphi } - \int_0^\tau \left( \intQ{ c_v \vr \varphi } \right) \circ 
\D W 
\end{equation}
holds 
for any $\tau \geq 0$ and any $\varphi \in C^1(\overline{Q})$, $\varphi \geq 0$.
\end{itemize}
\end{Definition} 

It is worth noting that the solutions introduced above are weak in the PDE sense - partial derivatives are interpreted in the sense of distributions - but strong in the stochastic sense - stochastic integral is considered on the original probability space. Our goal is to show the following result.  

\begin{Theorem} \label{TM1}

Let $Q \subset R^N$, $N=2,3$ be a bounded domain, 
\[
Q = \cup_{i = 1}^\infty \Ov{Q}_i, \ Q_i \ \mbox{domains},\ Q_i \cap Q_j = \emptyset \ \mbox{for}\ i \ne j.
\]
Suppose that $\vr_0, \ \vt_0 \in L^1(Q)$ are $\mathfrak{F}_0$-adapted random variables satisfying   
\[
0 < \underline{\vr} \leq \vr_0 \leq \Ov{\vr},\ 0 < \underline{\vt} \leq \vt_0 \leq \Ov{\vt} \ \prst-\mbox{a.s.}
\]
for some deterministic constants $\underline{\vr}$, $\Ov{\vr}$, $\underline{\vt}$, $\Ov{\vt}$ and such that 
\[
\vr_0|_{Q_i} = \vr_{0,i},\ \vt_0|_{Q_i} = \vt_{0,i}\ i = 1,2, \dots,\ \mbox{where}\ \vr_{0,i}, \ \vt_{0,i} \ \mbox{are constant.}
\] 

Then there exists a deterministic constant $\mathcal{E}_0$ such that for any $\mathcal{E} > \mathcal{E}_0$ there exists an 
$\mathfrak{F}_0-$adapted random field $\vm_0 \in L^\infty(Q; R^N)$,  such that 
\[
\intO{ \left[ \frac{1}{2} \frac{|\vm_0|^2}{\vr_0} + c_v \vr_0 \vt_0 \right] } \geq \mathcal{E} \ \prst-\mbox{a.s.},
\]
and the problem \eqref{I1}, \eqref{I4}, \eqref{M1} admits infinitely many weak solutions in $(0, \infty) \times Q$ with the initial data 
\[
\vr_0,\ \vm_0, \ E_0 = \frac{1}{2} \frac{|\vm_0|^2}{\vr_0} + c_v \vr_0 \vt_0.
\]

\end{Theorem}

The rest of the paper is devoted to the proof of Theorem \ref{TM1}.

\section{Proof of Theorem \ref{TM1}}
\subsection{Constant initial data}
\label{H}

We first assume that $\vr_0$, $\vt_0$ are positive (random) constants
admitting deterministic lower and upper bounds as in Theorem \ref{TM1}. Later we extend the result to piecewise constant data by ``pasting'' solutions together. 

\subsubsection{Solenoidal fields}

We look for solenoidal momentum fields $\vm$ with vanishing normal trace, meaning 
\begin{equation} \label{H1}
\intQ{ \vm \cdot \Grad \varphi } = 0 
\ \mbox{for any}\ \varphi \in C^1(\Ov{Q}).
\end{equation}
If we then set $\vr (t, \cdot) = \vr_0$ for any $t \geq 0$, in particular, the equation of continuity \eqref{M2} is automatically satisfied
and $ \intQ{ \vr (t, \cdot) }$ is a semimartingale.

\subsubsection{Temperature field}

Writing the internal energy equation as
\[
\D (\vr e) + \Div (e \vm) \dt = - (\vr e) \circ \D W - p \Div \left( \frac{\vm}{\vr} \right) \dt
\]
we realize that, since $\Div \vm = 0$ and $\vr = \vr_0$ constant, then the unique solution is given by 
\begin{equation} \label{H2}
e = \vr \vt,\ \vt = \vt_0 \exp( - W(t) ),
\end{equation}
where $\vt_0$ is the constant initial temperature. Obviously both $\vr$ and $\vt$ are continuous $\{ \mathfrak{F}_t \}_{t \geq 0}$-adapted 
semimartingales. When we have established that $\intQ{ |\vm|^2 \varphi }$ is a semimartingale, then also 
$\intQ{ E(t, \cdot) \varphi }$ is a semimartingale.

\subsubsection{Momentum equation}\label{Mo}

The computations of this subsection have to computed rigouroulsy in the reversed order, using the rules of Stratonovich calculus and starting from 
the $\mathfrak{F}_0 $-measurable process $\vc{v}(t,\cdot)$.
We first present them in this order for convenience of intuition.

In view of \eqref{H2}, the momentum equation reads
\[
\D \vm + \frac{1}{2} \vm \circ \D W + \Div \left( \frac{\vm \otimes \vm}{\vr_0} \right) \dt + 
\Grad \left( \vr_0 \vt_0 \exp \Big( - W(t) \Big) \right) \dt = 0. 
\]
Using the chain rule for the Stratonovich integral, we obtain 
\[
\begin{split}
\D \left[ \vm \exp \left( \frac{1}{2} W(t) \right) \right] &+ \exp \left( \frac{1}{2} W(t) \right)\Div \left( \frac{\vm \otimes \vm}{\vr_0} \right) \dt \\ &+ 
\Grad \left( \vr_0 \vt_0 \exp \left( - \frac{1}{2} W(t) \right) \right) \dt = 0. 
\end{split}
\]
In order to apply the convex integration argument, we need to recast the equation in a suitable way. Thanks to the computations above based on 
Stratonovih calculus, it easy to observe that by introducing a new variable $\vc{w}$, 
\[
\vc{w} = \vm \exp \left( \frac{1}{2} W(t) \right),
\]
we obtain the following PDE with random coefficients 
\[
\partial_t \vc{w} + \exp \left( - \frac{1}{2} W(t) \right)\left[ \Div \left( \frac{\vc{w} \otimes \vc{w}} {\vr_0} \right) + 
\Grad \left( \vr_0 \vt_0  \right) \right] = 0. 
\]
Moreover, introducing a new time variable 
\[
t \approx \int_0^t \exp \left( - \frac{1}{2} W(s) \right)\ {\rm d}s 
\]
we obtain the system 
\begin{equation} \label{H3}
\partial_t \vc{v} + \Div \left( \frac{\vc{v} \otimes \vc{v} }{\vr_0} + \vr_0 \vt_0 \mathbb{I} \right) = 0,\ 
\Div \vc{v} = 0, \ \vc{v}(0, \cdot) = \vc{v}_0,
\end{equation} 
which can be now treated in the ``convex integration framework''.
To allow  the ``pasting''  of piecewise constant solutions, the problem \eqref{H3} should be supplemented with ``do nothing'' boundary conditions,
specifically, its weak formulation reads:
\begin{equation} \label{H4}
\int_0^\infty \intQ{ \left[ \vc{v} \cdot \partial_t \bfphi + 
 \left( \frac{\vc{v} \otimes \vc{v} }{\vr_0} + \vr_0 \vt_0 \mathbb{I} \right) 
: \Grad \bfphi \right] } \dt = - \intO{\vc{v}_0 \cdot \bfphi(0, \cdot) }
\end{equation}
for any $\bfphi \in C^1_c([0, \infty) \times \Ov{Q}; R^N))$. As shown in the forthcoming section, problem \eqref{H4} admits infinitely many solutions for suitable initial data provided by the method of convex integration. 

\subsection{Convex integration}
\label{C}

To finally apply the method of convex integration, we reformulate the problem \eqref{H4}. Specifically, we replace \eqref{H4} by
\begin{equation} \label{C1}
\int_0^\infty \intQ{ \left[ \vc{v} \cdot \partial_t \bfphi + 
 \left( \frac{\vc{v} \otimes \vc{v} }{\vr_0} - \frac{1}{N} \frac{|\vc{v}|^2}{\vr_0}  \mathbb{I} \right) 
: \Grad \bfphi \right] } \dt = - \intQ{\vc{v}_0 \cdot \bfphi(0, \cdot) }
\end{equation}
for any $\bfphi \in \DC([0, \infty) \times R^N; R^N))$. In addition, we prescribe the energy 
\[
\frac{1}{2} \frac{|\vc{v}|^2}{\vr_0} = K_0,
\]
where $K_0 > 0$ is a positive random variable adapted to $\mathfrak{F}_0$. 

If $\vr_0$, $K_0$ were deterministic quantities, the nowadays standard method developed by De Lellis 
and Sz\' ekelyhidi in \cite{DelSze} would yield the existence of an initial datum $\vc{v}_0 \in L^\infty(\Omega; R^N)$ such that:
\begin{itemize}
\item 
\[
\Div \vc{v}_0 = 0,\ \vc{v}_0 \cdot \vc{n}|_{\partial \Omega} = 0;
\]
\item the problem \eqref{C1} admits infinitely many solutions $\vc{v}$ satisfying 
\begin{equation} \label{C2}
\int_0^\infty  \intQ{ \vc{v} \cdot \Grad \bfphi }\dt = 0 
\end{equation}
for all $\bfphi \in C^1_c([0, \infty) \times \Ov{Q}; R^N)$;
\item
\begin{equation} \label{C3}
\frac{1}{2} \frac{|\vc{v}_0 |^2}{\vr_0} = 
\frac{1}{2} \frac{|\vc{v}(t, \cdot)|^2}{\vr_0} = K_0 \ \mbox{a.a. in}\ \Omega \ \mbox{for any}\ t \geq 0.
\end{equation}

\end{itemize}

This result has been extended to the random setting in \cite{BrFeHo2017B}. Indeed, if 
$\vr_0$, $\vt_0$, and $K_0$ are $\mathfrak{F}_0-$adapted random variables, the stochastic version 
of the oscillatory lemma proved in \cite[Lemma 5.7]{BrFeHo2017B} can be applied to deduce that the solutions 
$\vc{v}$, obtained through process described in \cite{DelSze3},
are $\mathcal{F}_0$ adapted as random variables considered in the space 
$C_{{\rm weak}}([0, \infty); L^2(Q; R^N))$. More specifically, the random variable 
\[
t \mapsto \intQ{ \vc{v}(t, \cdot) \cdot \bfphi } \in C[0, \infty) 
\ \mbox{is}\ \mathfrak{F}_0-\mbox{adapted} 
\]
for any $\bfphi \in C^1(\Ov{Q}; R^N)$. Note that the present situation is much simpler
than in \cite{BrFeHo2017B} as the $\sigma$-field $\mathcal{F}_0$ is independent of time. 

Finally, we fix $K_0$ in such a way that 
\begin{equation} \label{C4}
\frac{N}{2} K_0 = \Lambda_0 - \vr_0 \vt_0 > 0
\end{equation}
where $\Lambda_0$ is a (random) constant. Note that, in view of our hypotheses imposed on the data $\vr_0$, $\vt_0$, the quantity 
$\Lambda_0$ can be chosen in a deterministic way. 

Summarizing we obtain the following result.

\begin{Proposition} \label{CP1}

Let $Q \subset R^N$, $N=2,3$ be a bounded domain. Suppose that $\vr_0$, $\vt_0$, $\Lambda_0$ are given (real valued) random variables that are 
$\mathfrak{F}_0$-adapted and satisfy 
\[
\vr_0,\ \vt_0, \ 
\Lambda_0 - \vr_0 \vt_0 > 0 \ \prst-\mbox{a.s.}
\]

Then there exists a random variable $\vc{v}_0 \in L^\infty_{{\rm weak(*)}}(Q; R^N)$ and infinitely many 
\[
\vc{v} \in C_{{\rm weak}}([0, \infty); L^2(Q; R^N))
\]
satisfying:

\begin{itemize}
\item
\begin{equation} \label{C5}
t \mapsto \intQ{ \vc{v}(t, \cdot) \cdot \bfphi } \in C[0, \infty) \ \mbox{are}\ \mathfrak{F}_0-\mbox{adapted} 
\end{equation}
for any $\bfphi \in C^1(\Ov{Q}; R^N)$;
\item 
\begin{equation} \label{C6}
\int_0^\infty  \intQ{ \vc{v} \cdot \Grad \bfphi }\dt = 0 
\end{equation}
for all $\bfphi \in C^1_c([0, \infty) \times \Ov{Q}; R^N)$;
\item
\begin{equation} \label{C7}
\frac{1}{2} \frac{|\vc{v}_0 |^2}{\vr_0} = 
\frac{1}{2} \frac{|\vc{v}(t, \cdot)|^2}{\vr_0} = \frac{2}{N} \left( \Lambda_0 - \vr_0 \vt_0 \right)\ \mbox{a.a. in}\ Q \ \mbox{for any}\ t \geq 0;
\end{equation}
\item 
\begin{equation} \label{C8}
\int_0^\infty \intQ{ \left[ \vc{v} \cdot \partial_t \bfphi + 
 \left( \frac{\vc{v} \otimes \vc{v} }{\vr_0} + \Big( \vr_0 \vt_0 - \Lambda_0 \Big)  \mathbb{I} \right) 
: \Grad \bfphi \right] } \dt = - \intQ{\vc{v}_0 \cdot \bfphi(0, \cdot) }
\end{equation}
for any $\bfphi \in C^1_c([0, \infty) \times \Ov{Q}; R^N))$. 

\end{itemize}

\end{Proposition}

\begin{Remark}

Note that the original result of De Lellis and Sz\' ekelyhidi \cite{DelSze3} would apply without modification should the initial data $\vr_0$, $\vt_0$ be deterministic.

\end{Remark}

The conclusion of Proposition \ref{CP1} should be seen as a starting point of the existence of infinitely many solutions claimed in Theorem \ref{TM1}. Note that, at this level, the density $\vr$, as well as the kinetic energy $\frac{|\vm|^2}{\vr}$ are in fact 
constants independent of the time variable.

\subsection{Piecewise constant data}
\label{P}

We claim that the conclusion of Proposition \ref{CP1} remains valid if the quantities $\vr_0$, $\vt_0$, and $\Lambda_0$ are piecewise constant random variables as required in Theorem \ref{TM1}. Specifically, we suppose that 
\[
\Ov{Q} = \cup_{i = 1}^\infty \Ov{Q}_i,\ Q_i \cap Q_j = \emptyset \ \mbox{if}\ 
i \ne j, 
\]
and 
\[
\vr_0 = \vr_{0,i},\ \vt_0 = \vt_{0,i},\ \Lambda_0 = \Lambda_{0,i} 
\ \mbox{in}\ Q_i,\ i = 1,\dots.
\] 
Indeed such a generalization is possible as the integrals in \eqref{C6}, \eqref{C8} are additive, and the test functions need not vanish on $\partial Q_i$. We simply apply Proposition \ref{CP1} on each $Q_i$ and take the sum of the corresponding integrals in 
\eqref{C6}, \eqref{C8}.

In addition, if $\vr_0$, $\vt_0$ are bounded by deterministic constants as in Theorem \ref{TM1}, we can choose the constants $\Lambda_{0,i} = \Lambda_0$ the same on each $Q_i$. In particular, equation \eqref{C8} gives rise to 
\begin{equation} \label{P1}
\int_0^\infty \intQ{ \left[ \vc{v} \cdot \partial_t \bfphi + 
\frac{\vc{v} \otimes \vc{v} }{\vr_0} : \Grad \bfphi + \vr_0 \vt_0 \Div \bfphi 
\right] } \dt = - \intQ{\vc{v}_0 \cdot \bfphi(0, \cdot) }
\end{equation}
for any $\bfphi \in \DC([0, \infty) \times \Ov{Q}; R^N))$ as long as 
 $\bfphi \cdot \vc{n}|_{\partial Q} = 0$.

Next, we derive from \eqref{C2} that  
\begin{equation} \label{P2}
\int_0^\infty \int_{Q} \vc{v} \cdot \Grad \varphi \ \dx \dt = 0 
\end{equation}
for any $\varphi \in C^1_c([0, \infty) \times \Ov{Q})$. 

Finally, the kinetic energy is piecewise constant and independent of time,  
\begin{equation} \label{P3}
\frac{1}{2} \frac{|\vc{v}(t, \cdot)|^2}{\vr} = \frac{1}{2} \frac{|\vc{v}_0|^2}{\vr_0} = \frac{2}{N} (\Lambda_0 - \vr_0 \vt_0) 
\ \mbox{a.a. in}\ Q \ \mbox{for any}\ t \geq 0.
\end{equation}

Thus Proposition \ref{CP1} can be extended to piece--wise constant data as follows. 

\begin{Proposition} \label{CP2}

Let $Q \subset R^N$, $N=2,3$ be a bounded domain, 
\[
\Ov{Q} = \cup_{i = 1}^\infty \Ov{Q}_i,\ Q_i \cap Q_j = \emptyset \ \mbox{if}\ 
i \ne j. 
\]
Suppose that $\vr_0,\ \vt_0 \in L^1(Q)$,  $\Lambda_0 \in R$ are given random variables that are 
$\mathfrak{F}_0$-adapted, with $\vr_0$, $\vt_0$ piecewise constant, meaning
\[
\vr_0 = \vr_{0,i} > 0,\ \vt_0 = \vt_{0,i} > 0 \ \mbox{in}\  Q_i,\  
\Lambda_0 - \vr_0 \vt_0 > 0 \ \prst-\mbox{a.s.}
\]

Then there exists an $\mathfrak{F}_0-$measurable random variable $\vc{v}_0 \in L^\infty_{{\rm weak(*)}}(Q; R^N)$ and infinitely many 
\[
\vc{v} \in C_{{\rm weak}}([0, \infty); L^2(Q; R^N))
\]
satisfying:

\begin{itemize}
\item
\begin{equation} \label{CC5}
t \mapsto \intQ{ \vc{v}(t, \cdot) \cdot \bfphi } \in C[0, \infty) \ \mbox{are}\ \mathfrak{F}_0-\mbox{adapted} 
\end{equation}
for any $\bfphi \in C^1(\Ov{Q}; R^N)$;
\item 
\begin{equation} \label{CC6}
\int_0^\infty  \intQ{ \vc{v} \cdot \Grad \bfphi }\dt = 0, \ \int_0^\infty \int_{Q_i} \vc{v} \cdot \Grad \bfphi \ \dxdt = 0,
\ i=1,\dots 
\end{equation}
for all $\bfphi \in C^1_c([0, \infty) \times \Ov{Q}; R^N)$;
\item
\begin{equation} \label{CC7}
\frac{1}{2} \frac{|\vc{v}_0 |^2}{\vr_0} = 
\frac{1}{2} \frac{|\vc{v}(t, \cdot)|^2}{\vr_0} = \frac{2}{N} \left( \Lambda_0 - \vr_0 \vt_0 \right)\ \mbox{a.a. in}\ Q \ \mbox{for any}\ t \geq 0;
\end{equation}
\item 
\begin{equation} \label{CC8}
\int_0^\infty \intQ{ \left[ \vc{v} \cdot \partial_t \bfphi + 
 \left( \frac{\vc{v} \otimes \vc{v} }{\vr_0} + \vr_0 \vt_0   \mathbb{I} \right) 
: \Grad \bfphi \right] } \dt = - \intQ{\vc{v}_0 \cdot \bfphi(0, \cdot) }
\end{equation}
for any $\bfphi \in C^1_c([0, \infty) \times \Ov{Q}; R^N))$, $\bfphi \cdot \vc{n}|_{\partial Q} = 0$. 

\end{itemize}

\end{Proposition}

\subsection{Rescaling time} 
\label{R}
In this last part of the proof, we show how to go back from the ``convex integration constructed $\vc{v}$'' to solutions $\vc{m}$ of the original system 
\eqref{I1}.
This can be done by computing and justifying formally the reversed transformations of the ones performes in Section \ref{Mo} 
to obtain the system for $\vc{v}$.
As a first step, we need to rescale time.
Consider the function of time, 
\[
\left< \vc{v}, \phi \right> \equiv 
t \mapsto \intQ{ \vc{v} (t, \cdot) \cdot \phi },\ t \in [0, \infty),\ \phi \in C^1(\Ov{Q}, R^N),\ \phi \cdot \vc{n}|_{\partial Q} = 0.
\]
It follows from \eqref{CC8} that $\left< \vc{v}, \phi \right>$ is globally Lipschitz on $[0, \infty)$ with the time derivative 
\[
\frac{{\rm d}}{\dt} \left< \vc{v}, \phi \right> =  \intQ{ \left[ \frac{\vc{v} \otimes \vc{v} }{\vr_0} : \Grad \phi + \vr_0 \vt_0 \Div \phi 
\right] } \ \mbox{for a.a.}\ t \in (0, \infty).
\]

We introduce a new function $\vc{w} \in C_{{\rm weak}}([0, \infty); L^2(\Omega; R^N))$, 
\[
\begin{split}
\left< \vc{w}, \phi \right> &\equiv 
\intQ{ \vc{w} (t, \cdot) \cdot \phi } = 
\intQ{ \vc{v} \left( \int_0^t \exp \left( - \frac{1}{2} W(s) \right) \D s, x \right) \phi (x) },\\ 
\phi &\in C^1(\Ov{Q}; R^N), \ \phi \cdot \vc{n}|_{\partial Q} = 0.
\end{split}
\]
Note carefully that $\vc{w}$ is $(\mathfrak{F})_{t \geq 0}$-adapted for any $\phi$, where $(\mathfrak{F})_{t \geq }$ is the filtration associated to the 
noise $W$.

Since $\left< \vc{v}; \varphi \right>$ is Lipschitz function of time, we can use the abstract chain rule (see e.g. Ziemer \cite{ZIE}) 
to deduce that, $\prst-$a.s.,
\begin{equation} \label{R1}
\frac{\D}{\dt} \intQ{ \vc{w} \cdot \bfphi } = \exp \left( - \frac{1}{2} W(t) \right)
\intQ{ \left[ \frac{\vc{w} \otimes \vc{w} }{\vr_0} : \Grad \bfphi + \vr_0 \vt_0 \Div \bfphi 
\right] } 
\end{equation}
for any $\bfphi \in C^1(\Ov{Q}; R^N), \ \bfphi \cdot \vc{n}|_{\partial Q} = 0$, $\vc{w}(0) \equiv \vc{w}_0 = \vc{v}_0$.

Finally, we observe that the relations \eqref{P2}, \eqref{P3} remain valid for $\vc{w}$,
specifically, 
\begin{equation} \label{R2}
\int_0^\infty \int_{Q_i} \vc{w} \cdot \Grad \varphi \ \dx \dt = 0,\ i=1, \dots  
\end{equation}
for any $\varphi \in C^1_c([0, \infty) \times \Ov{Q})$, and   
\begin{equation} \label{R3}
\frac{1}{2} \frac{|\vc{w}(t, \cdot)|^2}{\vr} = \frac{1}{2} \frac{|\vc{w}_0|^2}{\vr_0} = \frac{2}{N} (\Lambda_0 - \vr_0 \vt_0) 
\ \mbox{a.a. in}\ Q \ \mbox{for any}\ t \geq 0.
\end{equation}

\subsection{Chain rule for Stratonovich integral}

\label{S}
Finally, we can introduce the momentum $\vc{m}$ in terms of the rescaled $\vc{v}$ (i.e. of $\vc{w}$), so to get solutions to system 
\eqref{I1}. To this aim it is essential the use of Stratonovich calculus.
We introduce the momentum 
\[
\vc{m} = \vc{w} \exp\left( - \frac{1}{2} W(t) \right)
\]
noting that 
\[
\vc{m}(0, \cdot) \equiv \vc{m}_0 = \vc{w}_0 = \vc{v}_0.
\]
Obviously, the relation \eqref{R2} applies to $\vc{m}$, 
\begin{equation} \label{S1}
\int_0^\infty \int_{Q_i} \vc{m} \cdot \Grad \varphi \ \dx \dt = 0,\ i=1, \dots  
\end{equation}
for any $\varphi \in C^1_c([0, \infty) \times \Ov{Q})$. 

Using the basic properties of Stratonovich integral, we obtain 
\begin{equation} \label{S2}
\D \left( b \exp (W) \right) = \exp(W) \D b + b \exp(W) \circ \D W
\end{equation}
whenever $b$ is a Lipschitz function.

At this stage, we are ready to finish the proof of Theorem \ref{TM1}.

\subsubsection{Equation of continuity}

Setting $\vr = \vr_0$ and using \eqref{S1} we easily deduce the equation of continuity 
\begin{equation} \label{S3}
\D \intQ{ \vr \varphi } = \intQ{ \vc{m} \cdot \Grad \varphi } \dt,\ \intO{ \vr(0, \cdot) \varphi }
= \intQ{ \vr_0 \varphi } 
\end{equation}
for any $\varphi \in C^1(\Ov{Q})$.

\subsubsection{Internal energy, entropy, total energy}

We define 
\[
\vt = \vt_0 \exp(- W(t))
\]
and, using the relation \eqref{S1}, \eqref{S2}, we easily deduce the internal energy equation 
\begin{equation} \label{S4}
\begin{split}
\D \intQ{ c_v \vr \vt \varphi } &= \intQ{ c_v \vr \vt \frac{\vc{m}}{\vr} \cdot \Grad \varphi } \dt
- \intQ{ c_v \vr \vt \varphi } \circ \D W\\
\intQ{ \vr \vt (0, \cdot) \varphi } &= \intQ{ \vr_0 \vt_0 \varphi }
\end{split}
\end{equation}
for any $\varphi \in C^1(\Ov{Q})$.

Similarly, seeing that the entropy is 
\[
s(\vr, \vt) = c_v \log (\vt) - \log(\vr) = c_v \log(\vt_0) - \log(\vr_0) - c_v W, 
\]
we obtain the entropy equation 
\begin{equation} \label{S5}
\begin{split}
\D \intQ{ \vr s( \vr, \vt) \varphi } &= \intQ{ s( \vr ,\vt) {\vc{m}} \cdot \Grad \varphi } \dt
- \intQ{ c_v \vr \varphi } \circ \D W\\
\intQ{ \vr s(\vr, \vt) (0, \cdot) \varphi } &= \intQ{ \vr_0 s(\vr_0, \vt_0) \varphi }
\end{split}
\end{equation}
for any $\varphi \in C^1(\Ov{Q})$.

\begin{Remark}

Note that we have shown the existence of infinitely many solutions that satisfy the entropy \emph{equation} instead of the mere 
inequality required in Definition \ref{D1}.

\end{Remark}

Finally, by virtue of \eqref{R3}, the total energy reads  
\[
\begin{split}
E &= \frac{1}{2} \frac{|\vc{m}|^2}{\vr} + c_v \vr \vt = 
\exp (-W(t)) \left( \frac{1}{2} \frac{|\vc{w}|^2}{\vr} + c_v \vr_0 \vt_0 \right)\\ 
&= \exp (-W(t)) \left( \frac{2}{N} (\Lambda_0 - \vr_0 \vt_0 ) + c_v \vr_0 \vt_0 \right). 
\end{split}
\]
Thus, similarly to the above, we deduce the total energy balance 
\begin{equation} \label{S6}
\begin{split}
\D \intQ{ E \varphi } &= \intQ{ \left( E + \vr \vt \right) \frac{\vc{m}}{\vr} \cdot \Grad \varphi } \dt
- \intQ{ E \varphi } \circ \D W\\
\intO{ E (0, \cdot) \varphi } &= \intQ{ \left( \frac{1}{2} \frac{ |\vc{m}_0|^2}{\vr_0} + c_v \vr_0 \vt_0 \right) \varphi } 
\end{split}
\end{equation}
for any $\varphi \in C^1(\Ov{Q})$.

\subsubsection{Momentum equation}

We multiply \eqref{R1} on $\exp \left( - \frac{1}{2} W \right)$ obtaining
\[
\begin{split}
\exp \left( - \frac{1}{2} W(t) \right) \frac{\D}{\dt} \left[ \exp \left( \frac{1}{2} W(t) \right) \intQ{ \vc{m} \cdot \bfphi } \right] = 
\intQ{ \left[ \frac{\vc{m} \otimes \vc{m} }{\vr} : \Grad \phi + \vr \vt \Div \bfphi 
\right] },
\end{split}
\]
which, in view of \eqref{S2}, gives rise to the desired conclusion 
\begin{equation} \label{S7}
\begin{split}
\D \intO{ \vc{m} \cdot \bfphi }  = 
\intQ{ \left[ \frac{\vc{m} \otimes \vc{m} }{\vr} : \Grad \bfphi + \vr \vt \Div \bfphi 
\right] } \dt - \frac{1}{2} \intQ{ \vc{m} \cdot \bfphi } \circ \D W 
\end{split}
\end{equation}
for any $\bfphi \in C^1(\Ov{Q}; R^N), \ \bfphi \cdot \vc{n}|_{\partial Q} = 0$.

\bigskip

We have shown Theorem \ref{TM1}.
\section{Appendix}
Since the use of Stratonovich integrals and calculus is an essential tool of
this work, we summarize some definitions and rules;\ everything can be found
in details in Kunita \cite{Kuni}. Given a probability basis $\left(  \Omega
,\mathcal{F},P\right)  $ with a right-continuous complete filtration $\left(
\mathcal{F}_{t}\right)  _{t\geq0}$, the general concept of \textit{continuous
semimartingale} can be found in many books, see e.g. Kunita \cite{Kuni},
Revuz and Yor \cite{RevuzYor}. Examples of continuous semimartingales are the Brownian
motion $\left(  \beta_{t}\right)  _{t\geq0}$, the deterministic (Riemann type)
integrals $\int_{0}^{t}X_{s}ds$ of continuous semimartingales $\left(
X_{t}\right)  _{t\geq0}$ and three objects we now define. Given two continuous
semimartingales $\left(  X_{t}\right)  _{t\geq0}$ and $\left(  Y_{t}\right)
_{t\geq0}$, the following limits of finite Riemann type sums exist, understood
as limits in probability. Given $t>0$, let $\left(  \pi_{n}\right)
_{n\in\mathbb{N}}$ be a sequence of partitions of $\left[  0,t\right]  $ and
denote points of $\pi_{n}$ generically by $t_{i}$. When we write $\sum
_{t_{i}\in\pi_{n}}$ we understand that the sum is extended to all indexes that
are admissible in the expression (since also $t_{i+1}$ appears). Then we set:%
\begin{align*}
\int_{0}^{t}X_{s}dY_{t}  & :=\lim_{n\rightarrow\infty}\sum_{t_{i}\in\pi_{n}%
}X_{t_{i}}\left(  Y_{t_{i+1}}-Y_{t_{i}}\right)  \\
\int_{0}^{t}X_{s}\circ dY_{t}  & :=\lim_{n\rightarrow\infty}\sum_{t_{i}\in
\pi_{n}}\frac{X_{t_{i}}+X_{t_{i+1}}}{2}\left(  Y_{t_{i+1}}-Y_{t_{i}}\right)
\\
\left[  X,Y\right]  _{t}  & :=\lim_{n\rightarrow\infty}\sum_{t_{i}\in\pi_{n}%
}\left(  X_{t_{i+1}}-X_{t_{i}}\right)  \left(  Y_{t_{i+1}}-Y_{t_{i}}\right)  .
\end{align*}
We call them It\^{o} integral, Stratonovich integral and covariation,
respectively. These limit, always in probability, can be also understood
uniformly in time on finite intervals, with due modification of the notations,
so that the partitions are not adapted to a single interval $\left[
0,t\right]  $. The class of continuous semimartingales is closed also under
the previous three operations;\ and under sum, product and in general
composition by functions $f\in C^{1,2}\left(  \left[  0,T\right]
\times\mathbb{R}^{d}\right)  $. 

When $\left(  Y_{t}\right)  _{t\geq0}$ is a Brownian motion $\left(  \beta
_{t}\right)  _{t\geq0}$, the integral $\int_{0}^{t}X_{s}d\beta_{t}$ is the
classical It\^{o} integral. It exists also when $X$ is just a continuous
adapted process; and with an alternative definition it is well defined also in
more general cases; in the framework of It\^{o} calculus, functions $f$ of
class $C^{1,2}\left(  \left[  0,T\right]  \times\mathbb{R}^{d}\right)  $
suffice to write a chain rule. On the contrary, the two objects $\int_{0}%
^{t}X_{s}\circ dY_{t}$ and $\left[  X,Y\right]  _{t}$ are quite restrictive
from the viewpoint of existence and the class of continuous semimartingales
looks the right one for a general theory; and manipulations require functions
$f$ of class $C^{1,3}\left(  \left[  0,T\right]  \times\mathbb{R}^{d}\right)
$. This is the price to work with them. The advantage are the rules of
calculus. These rules (summarized by the multidimensional chain rule) based on
It\^{o} integrals are well known to be modified by the presence of a
correction term. When Stratonovich integral is used, the rules are the same as
deterministic calculus. For instance, in this work we use the fact that, for
two continuous semimartingales $\left(  X_{t}\right)  _{t\geq0}$ and $\left(
Y_{t}\right)  _{t\geq0}$,
\[
X_{t}Y_{t}=X_{0}Y_{0}+\int_{0}^{t}X_{s}\circ dY_{s}+\int_{0}^{t}Y_{s}\circ
dX_{s}%
\]
which is the rigorous formulation of the identity commonly written as%
\[
d\left(  X_{t}Y_{t}\right)  =X_{t}\circ dY_{t}+Y_{t}\circ dX_{t}.
\]
The analogous result with It\^{o} integrals is%
\[
X_{t}Y_{t}=X_{0}Y_{0}+\int_{0}^{t}X_{s}dY_{s}+\int_{0}^{t}Y_{s}dX_{s}+\left[
X,Y\right]  _{t}.
\]
More generally, if $f\in C^{1,2}\left(  \left[  0,T\right]  \times
\mathbb{R}^{d}\right)  $ and $X_{t}=\left(  X_{t}^{1},...,X_{t}^{n}\right)  $
is a vector of continuous semimartingales, then%
\[
df\left(  t,X_{t}\right)  =\left(  \partial_{t}f\right)  \left(
t,X_{t}\right)  dt+\sum_{i=1}^{n}\left(  \partial_{x_{i}}f\right)  \left(
t,X_{t}\right)  \circ dX_{t}^{i}%
\]
opposite to It\^{o} formula%
\begin{align*}
df\left(  t,X_{t}\right)    & =\left(  \partial_{t}f\right)  \left(
t,X_{t}\right)  dt+\sum_{i=1}^{n}\left(  \partial_{x_{i}}f\right)  \left(
t,X_{t}\right)  dX_{t}^{i}\\
& +\frac{1}{2}\sum_{i,j=1}^{n}\left(  \partial_{x_{i}}\partial_{x_{j}%
}f\right)  \left(  t,X_{t}\right)  d\left[  X^{i},X^{j}\right]  _{t}.
\end{align*}
A technical remark:\ when It\^{o} interpretation of integrals is given, the
integrand is just required to be continuous adapted hence $\left(
\partial_{x_{i}}f\right)  \left(  t,X_{t}\right)  $, $\left(  \partial_{x_{i}%
}\partial_{x_{j}}f\right)  \left(  t,X_{t}\right)  $ are admissible
integrands. When Stratonovich interpretation of integrals is chosen, the
integrand must be a continuous semimartingale. Hence $\left(  \partial_{x_{i}%
}f\right)  \left(  t,X_{t}\right)  $ must have such property and, by It\^{o}
calculus, we know it is true when $\partial_{x_{i}}f\in C^{1,2}\left(  \left[
0,T\right]  \times\mathbb{R}^{d}\right)  $. This is why the property $f\in
C^{1,3}\left(  \left[  0,T\right]  \times\mathbb{R}^{d}\right)  $ is required
in Stratonovich calculus.

\bigskip

\textbf{Acknowledgement.} 
The work has been essentially discussed during the stay of E.F. at Scuola Normale Superiore in Pisa, whose support and hospitality 
is gladly acknowledged.


\def\cprime{$'$} \def\ocirc#1{\ifmmode\setbox0=\hbox{$#1$}\dimen0=\ht0
  \advance\dimen0 by1pt\rlap{\hbox to\wd0{\hss\raise\dimen0
  \hbox{\hskip.2em$\scriptscriptstyle\circ$}\hss}}#1\else {\accent"17 #1}\fi}

\end{document}